\documentclass[12pt]{article}
\usepackage{amssymb,amsmath,amsfonts,latexsym}
\begin{document}
\title{On a variational approach to truncated problems of moments.}
\author{C.-G. Ambrozie\footnote{Supported by grants IAA100190903  GAAV and 201/09/0473 GACR, RVO: 67985840}, \mbox{ }Prague}
\date{}

\newtheorem{remark}{Remark}
\newtheorem{remarks}[remark]{Remarks}
\newtheorem{proposition}[remark]{Proposition}
\newtheorem{examples}[remark]{Examples}
\newtheorem{example}[remark]{Example}
\newtheorem{lemma}[remark]{Lemma} 
\newtheorem{corollary}[remark]{Corollary}
\newtheorem{theorem}[remark]{Theorem}

\maketitle
\begin{abstract}

We characterize the existence of the $L^1$ solutions of the truncated moments problem in several real variables on unbounded supports by the existence of the maximum of  certain concave Lagrangian functions.
% associated dually with an entropy optimization problem. 
A natural regularity assumption on the support is required.
 

\end{abstract}

\section{Introduction}

%One venue, preferred by physicists
%and statisticians, is to add to the moment constraints the maximum entropy assumption

The present paper is concerned with the truncated problem of moments in several real variables, in the following context.
Let $n\in \mathbb{N}$ and fix a closed subset $T\! \not =\! \emptyset$ of $ \mathbb{R}^n$, 
a finite subset $I\!\subset\! (\mathbb{Z}_{+})^n$  with $0\!\in\! I$ and a set $g\! =\! (g_i )_{i\in I}$  of real numbers with $g_0 \! =\! 1$,
 where $\mathbb{Z}_{+}\! =\! \mathbb{N} \cup  \{ 0\}$.
Typically  a problem of moments  \cite{1} requires to establish if there exist Borel measures $\nu\! \geq \! 0$ on $\mathbb{R}^n$, supported on $T$,
such that $\int_T |t^i |d\nu (t)<\infty$ and
$
\int_{T}t^i d\nu (t) =g_i 
$
 for all $i\in I$. As usual $t^i =t_{1}^{i_1}\cdots t_{n}^{i_n}$  where $t=(t_1 ,\ldots ,t_n )$ is the variable in $\mathbb{R}^n$ and $i=(i_1 ,\ldots ,i_n )$ is a multiindex. In this case we call $\nu$ a 
{representing measure} of $g$, and $g_i$  the {moments} of $\nu$. %In what follows
 We are interested  in those measures $\nu =fdt$ that are absolutely continuous with respect to the $n$-dimensional Lebesgue measure $dt=dt_1 \cdots dt_n$, in which  case we call $f$ a 
{representing density} of $g$. 
 Namely the (class of equivalence of the) Lebesgue integrable function $f$ is $\, \geq 0$ almost everywhere ({a.e.}) on $T$, has finite moments of orders $i\in I$ and
\begin{equation}
\label{prima}
 \int_T t^i f(t) dt =g_i \mbox{ }\mbox{ }\, (i\in I).
  \end{equation}
Given  partial information in the integral form $\int_T t^i \, {\rm f}\,  \rho dt =g_i$ about representing densities ${\rm f}$ on a probability space $(T,\rho dt)$, endowed with a  reference density  $\rho $, does not determine them uniquely.
 An approach favorite to  physicists
and statisticians
 is to choose that particular density ${\rm f}_*$, minimizing 
 the entropy functional $h({\rm f} )\! =\! \int_T ({\rm f}\ln {\rm f}) \rho \, dt$ amongst all solutions
 of the moments constraints. This uniquely selects the unbiased
 probability distribution
 ${\rm f}_* $ (that proves to have the form ${\rm f}_* (t)=e^{\sum_{i\in I}\lambda_{i}^* t^i }$)  on the  knowledge of the prescribed average values $g_i$ of  $t^i$, 
where $t$ is considered as a $T$-valued random variable 
with repartition $\rho$  \cite{12'}, \cite{6}, \cite{12''}, \cite{kull}. Under suitable hypotheses, ${\rm f}_*$ turns to exist, even for  measures
more general than $\rho dt$.
A main tool to this aim is Fenchel duality
 \cite{BL},    \cite{Moreau}, \cite{rod}, \cite{ro}, that deals with 
 minimizing  convex functions $h:X\!\to \!\mathbb{R}\!\cup\! \{ \infty\} $ on convex subsets of locally convex spaces $X$, in connection
 with the dual problem of maximizing $-h^*$, where $h^* :X^* \!\to\! \mathbb{R}\!\cup\! \{  \infty\}$ is the convex conjugate of $h$, called also 
its  Legendre-Fenchel transform  \cite{rod}, \cite{ro},  defined on the dual $X^*$ of $X$ by 
 $h^* (y)\! =\! \sup \{ \langle x,y\rangle\! -\! h (x)\! :\! h(x)\! <\! \infty \} $. Typically $\inf h=\max (-h^* )$ and, briefly speaking, 
minimizing $\int_T {\rm f}\ln {\rm f} \rho dt $ as above
 is  to
find $\lambda^* =(\lambda_{i}^* )_{i\in I}$ maximizing 
$
L(\lambda )=\sum_{i\in I}g_i \lambda_i -\int_T e^{\sum_{i\in I}\lambda_i t^i }\rho dt
$.
% of $\lambda =(\lambda_i )_{i\in I}$.
Many  results exist in this direction
 \cite{MaEn},   \cite{bla} -- \cite{6}, \cite{Levermore}, \cite{junk}, \cite{Leww} -- \cite{MP}.
Additional hypotheses are always  necessary when  the conclusion $\inf h=\min h$ is sought for,
since there are $g$ for which the primal attainment fails   \cite{Levermore}, \cite{junk} although  problem (\ref{prima}) has solutions.
 
By  Theorem \ref{marg} we prove  that the feasibility of problem (\ref{prima}) is equivalent to the boundedness from above $\sup L<\infty$
 with attainment $\sup L=\max L$ for the concave Lagrangian function $L$. This holds no matter whether $\inf h$ is attained or not
(the general theory still provides us with $\inf h=\max L$).

Initiated by Stieltjes, Hausdorff,  Hamburger and Riesz, the area of the truncated problems of moments  
nowadays knows various other approaches, 
based for instance on operator methods or sums-of-squares representations for positive polynomials   
 \cite{7} --    \cite{11},   \cite{KM},  \cite{14''}.  Although important, these  topics remain beyond the aim of this work, focused on our mentioned Theorem \ref{marg}.
 
The author got the idea to consider $L$ instead of $h$ from the works  \cite{bla} where a similar characterization exists, and \cite{Levermore}, \cite{junk}, drawn to
 his attention by professor Mihai Putinar. Our statement and proof are rather general,  independent of these cited works.
% They  rely
%on  [Corollary 2.6, \cite{BL}] for the dual attainment of $L$ (the implication (a) $\Rightarrow$ (b) in Theorem \ref{poz}) and on  [Theorem 7,\cite{RHav}] 
%for the feasibility of (\ref{prima}) (the opposite implication (b) $\Rightarrow$ (a)). 

\section{Main results}

%Let $T$, $I$, $g$ be as stated in the Introduction.
% Let $I\not =\emptyset$ be a finite subset of $\mathbb{Z}_{+}^n$. We call $I$ {\it regular} \cite{RHav} if be a closed set below some topics related to the general problem of moments. 
%Given    and  , a 
We remind  that a linear  Riesz functional  $\varphi_\gamma $ \cite{CuHr} associated to a set $\gamma =(\gamma_i )_{i\in J}$ of real numbers 
$\gamma_i $ for $J\subset \mathbb{Z}_{+}^n$ is defined on the polynomials $p$ from 
the linear span  of $X_{1}^{i_1}\ldots X_{n}^{i_n}$ where $i=(i_1 ,\ldots ,i_n )\in J$ by
 $\varphi_\gamma  X^i = \gamma_i$. One calls $\varphi_\gamma$  $T$-{positive}  \cite{CuHr} if $\varphi_\gamma p\geq 0$ whenever $p(t)\geq 0$ for all $t\in T$.
 If $\gamma$ has   representing measures $\nu \geq 0$ on $T$,  $\varphi_\gamma $  is $T$-positive since  $\varphi_\gamma p=\int_T p\, d\nu $ for any 
such polynomial $p$. % polynomial.
% in th $is case.
 In the full case $J=\mathbb{Z}_{+}^n$ 
the $T$-positivity condition is  sufficient for the existence of the representing measures, 
by the  Riesz-Haviland theorem \cite{h}. %, by the  Riesz-Haviland theorem \cite{h}
% is concerned with the full problem of moments in $n$ real variables when all moments $g_i$ are prescribed, for $I=\mathbb{Z}_{+}^n$.
%As usual, $t^i =t_{1}^{i_1}\cdots t_{n}^{i_n}$ and $|i|=i_1 +\cdots +i_n$ for $i=(i_1 ,\ldots ,i_n )$
An analogue of this theorem  \cite{CuHr} for the truncated case   $I\! =\!\{ i\! :\! |i|\! \leq \! 2k\} $ 
characterizes  the existence of the representing measures by
the existence of 
$T$-positive extensions  of $\varphi_\gamma$ to the space of polynomials of degree $\leq 2k+2$. % \cite{CuHr}.
%, namely extensions $\tilde{\varphi}_{g}$ such that $\tilde{\varphi}_g p\geq 0$ for all $p$ with $\mbox{deg}\, p\leq 2k+2$ such that $p\geq_{\, T}0$.
For later use, we state below  a version  of these results  (Theorem \ref{poz}) and a Fenchel theoretic result of dual attainment 
  (Theorem \ref{maxim}).
%from below, respectively.
% (see \cite{h} for the full case). 
\vspace{1 mm}

\noindent {\bf Definitions} We call $T\!$ {\it regular} \cite{RHav} if
for any $t\in T$ and $\varepsilon >0$ the Lebesgue measure of the set $\{ x\in T: \| x-t\| <\varepsilon \}$ is positive. 
As usual $\| t \| =(\sum_{\iota =1}^n t_{\iota}^2 )^{1/2}$.  For any $i\in I$ set $\sigma_i = \{  j \in \mathbb{Z}_{+}^n  : j_k \! =\! \mbox{either } 0 \mbox{ or } i_k ,\, 1\leq k\leq n\}$. We call $I$ {\it regular} \cite{RHav} if
$\sigma_i \subset I$ for all $i\in I$. 
Define  $\Gamma ,\, G \subset \mathbb{R}^N$ ($N=\mbox{card}\, I$) by 
  $\Gamma \! =\! \{ \gamma \! =\! (\gamma_i )_{i\in I}:\,  \exists \mbox{ } \mbox{\rm  measures } \nu \geq 0\, \mbox {\rm on}\, T\, \mbox{\rm with}\,
  \int_T t^i d\nu (t)\! =\! \gamma_i ,\, \, i\in I\} $
  and
  $
  G\! =\! \{ \gamma\! =\! (\gamma_i )_{i\in I} \not  =\! 0:\, \exists \mbox{ }\, f\in L_{+}^1 (T,dt)\, \mbox{ }\mbox{\rm such that}\,
  \int_T t^i f(t)dt\! =\! \gamma_i ,\, i\in I\} .
  $
% with $|t^i | $  assumed to be in $L^1$. 
 The notation $L^p (T,\mu )$,  $L^p (\mu )$ for $\mu$  measure on $T$, $1\!\leq \! p\! \leq\! \infty$ has the usual meaning. In particular $L_{+}^1 (T, \mu )$  is the set of all 
$f\in L^1 (T,\mu )$, $f\! \geq\! 0$ 
$\mu$-a.e. % on $T$. 
For $\gamma \! =\! (\gamma_i )_{i\in I}$, $\varphi_\gamma $ is 
the linear functional defined on the  span $P_I $ $\!\subset \!\mathbb{R}[X_1 ,\ldots ,X_n ]$ of all  $X^i$ with $i\!\in \! I$  
   by
 $\varphi_\gamma  X^i \! =\!  \gamma_i$.
Set   ${\rm e}_\iota \! =\! (0,\ldots ,\stackrel{\iota}{1},\ldots ,0)$ for $1\! \leq \!\iota\!\leq\! n$.
 
By [Theorem 6,\cite{RHav}] the convex cone $G$  is the dense interior of the  cone $\Gamma$. % of all data $g$ that have representing measures.
 
 \begin{theorem}
 \label{poz} {\em [Theorem 7,\cite{RHav}]}
 Let $T\subset \mathbb{R}^n$ be a closed regular set, $I\subset \mathbb{Z}_{+}^n$  a finite regular set and $g=(g_i )_{i\in I}$  a set of numbers with $g_0 =1$. Then 
$g\in G$ $\Leftrightarrow$ 
  $\varphi_g \, p>0$ for every  $p\in P_I \setminus \{ 0\}$ such that $p(t)\geq 0$ for all $t\in T$.
 \end{theorem} 
 
%Whenever used, the norm on $\mathbb{R}^n$ or $\mathbb{R}^N$
 %($N=\mbox{card}\, I$) is the Euclidian one.

 \begin{theorem}
\label{maxim}
{\em [Corollary 2.6,\cite{BL}]}
Let $\mathcal{T}$ be a space with finite measure $\mu \! \geq \! 0$, $1\! \leq\!  p\! \leq\! \infty$ and $a_i \! \in \! L^{q} (\mu ) $, $g_i \! \in \!
\mathbb{R}$ for $i\! \in\! I$ 
$\! =\, $finite where $\frac{1}{p}\! +\! \frac{1}{q}\! =\! 1$. Let  $\phi :\mathbb{R}\! \to\!  (-\infty ,\infty ]$ 
be proper, convex, lower semicontinuous with $\phi |_{(0,\infty )}\! <\! \infty $.
If there are $x\! \in\!  L^p (\mu  )$, $x\! >\! 0$ a.e. such that $\phi \! \circ \! x\! \in\!  L^1 (\mu  )$ and
$ \int_\mathcal{T} a_i \, x \, d\mu \! =\! g_i $, then the quantities 
%$P\in [-\infty , \infty )$ and $D\in [-\infty , \infty \, ]$ defined respectively by
$$
P\! =\!  \inf \{ \! \int_\mathcal{T} \!  \phi (x(t)) d\mu (t)  :  x\in\! L^p (\mu ),   x\geq 0  \, \mbox{a.e.},\, 
\phi \!\circ\! x\in\! L^1 (\mu  ), \! \int_\mathcal{T} \! \!  a_i xd\mu \! =\! g_i \,  \forall i \, \} ,
$$
$$
D\! =\! \max  \{ \sum_{i\in I} g_i \lambda_i -  \int_\mathcal{T}  \! \phi^*   (   \sum_{i\in I} \lambda_i a_i (t)) \, d\mu (t) :
  \lambda_i  \in \mathbb{R}, \, \phi^* \! \! \circ\!  \sum_{i\in I} \lambda_i a_i \in L^1 (\mu ) \, \}
$$
are equal, $-\infty \leq P=D<\infty $ and the maximum  $D$ is attained.  
\end{theorem}
 
%\section{Main result}

Theorem \ref{marg} is a reminiscent to  [Theorem 4, \cite{MaEn}], where $\int_T {\rm f}\ln {\rm f} \rho dt$ is minimized subject to $\int_T t^i {\rm f} \rho dt =g_i$ under stronger hypotheses on $\rho$, like  $\rho (t)\sim e^{-\varepsilon \| t\|^p}$ with $p>2k$ (to fit the notation in \cite{MaEn}, let $a=1$ and  our $f:=\rho \, {\rm f}$, whence $L_{\rho, a, g}(\lambda )=L(\lambda -\lambda_0 )+1$, with $\lambda_0 =(\lambda_{0 i})_{i\in I}$ where $\lambda_{0i}=\delta_{i,0}$
and $\delta_{i,j}$ is Kronecker's symbol,
$\delta_{i,j}\! =\! 1$ if $i\! =\! j$ and $0$ if $i\! \not\! =\! j$).  
%and with partly different proof.
%, characterizes the existence of the solutions $f$ of (\ref{prima}). 
Although we do not obtain here the existence of a maximum entropy solution $f_*$, our present hypothesys on $\rho$ are weaker, while  condition $g\in G$  still characterized  in Lagrangian terms.  Our proof below relies on Theorem \ref{poz} ([Theorem 7,\cite{RHav}]) and Theorem \ref{maxim} ([Corollary 2.6,\cite{BL}]).

%For  various  topics  on convex functions to be used in what follow, we refer also to Rockafellar's book  \cite{ro}. 
%This  makes it interesting even if the conclusions we obtain below are weaker, too.

 \begin{theorem}
 \label{marg}
 Let $T\subset \mathbb{R}^n$ be a closed regular set. Let
 $I\subset \mathbb{Z}_{+}^n$ be a finite regular set  such that $\max_{i\in I}|i|\! =\! 2k $  where  $k\! \in\!  \mathbb{N}$. Assume    $\, 2k{\rm e}_\iota \! \in\! I$ $(1\! \leq \! \iota\! \leq \! n)$.  Let   
%, containing also $n$ multiindices of the form $(0,\ldots ,\stackrel{j}{2k_j},\ldots ,0)$ with $k_j \geq 1$ for $j=\overline{1,n}$.
$g=(g_i )_{i\in I}$  be a set of numbers with $g_0 =1$. 
Fix $\rho \in L^1 (T,dt)$, $\rho >0$ a.e. % such that $\int_T \| t\|^{\max_{i\in I}|i|} \rho (t)dt<\infty$.
 The following statements {\em (a)} and {\em (b)} are equivalent:
 
{\em (a)} There exist functions  $f\in L_{+}^1 (T,dt)$ such that $\int_T |t^i |f(t)dt<\infty$ and
$$
 \int_T t^i f(t)dt=g_i \mbox{ }\mbox{ }\mbox{ }(i\in I);
$$
% \end{equation}

{\em (b)} The  functional $L :\mathbb{R}^N \to  \mathbb{R}\cup \{ -\infty \}$ defined by
$$L(\lambda )=\sum_{i\in I}g_i \lambda_i -\int_T e^{\sum_{i\in I}\lambda_i t^i }\rho (t)dt,\mbox{ }\mbox{ }\mbox{ }\mbox{ }\, 
\lambda =(\lambda_i )_{i\in I} $$
  is bounded from above and  $\, \sup L$  is attained in a (unique) point $\lambda^* $.

 \end{theorem}

   {\em Proof}.    Since $L(0)>-\infty$, $L\not \equiv -\infty$. Since $g_0 =1$, each of the conditions (a) and (b) implies that  $T$ has positive Lebesgue measure, finite or not.
   Hence by means of Jensen's inequality one can show that $L$ is strictly concave. Then whenever $\sup L$ is finite and attained at some point $\lambda^*$, 
this $\lambda^*$ is unique.

(a) $\Rightarrow$ (b) %We  use  Fenchel theoretic results  \cite{BL},  \cite{junk}, \cite{Levermore}.
The regularity condition on $T$ is not necessary for this implication. 
Let $\mu =\tilde{\rho} dt$ be the  measure  on $T$ with density $\tilde{\rho}:=\rho e^{-\sum_{\iota =1}^n t_{\iota}^{2k}}$. Then $0<\mu (T)<\infty$.
 Since (\ref{prima}) has a  solution $f$, then $\tilde{f}:=f/\tilde{\rho}$ satisfies
\begin{equation}
\label{cheaptrick}
\int_T t^i \tilde{f} (t)\,  d\mu (t)=g_i \,  \mbox{  }\, \mbox{ (}i\in I\mbox{)}.
\end{equation}
By [Theorem 2.9, \cite{BL}], see also [Lemma 4, \cite{RHav}] for $\beta =0$,  problem (\ref{cheaptrick})
has also a solution $f_0 \in  L^\infty (T)$ with $f_0 >0$ a.e. The conclusion $\sup L<\infty$ may hold either directly by Theorem \ref{maxim}, or by 
an elementary argument as shown below.
 Let $x=f_0 (t)$ a.e. and $y=\| f_0 \|_\infty +1 $
in the inequalities  $-e^{-1}\leq x\ln x\leq y\ln y$ for $0\leq x\leq y$, $y\geq 1$, then 
 integrate with respect to $\mu$. Hence  $f_0 \ln f_0  \in L^1 (T, \mu )$.
  Fix $\lambda =(\lambda_i )_{i\in I}$. 
 Let  $x=f_0 (t)$ and $y=\sum_{i\in I}\lambda_i t^i $ in the 
simple version $x\ln x -x \geq xy-e^{y}$ of Fenchel's inequality \cite{ro}, then  integrate. It follows, 
using  (\ref{cheaptrick}) for $f_0$,
that $$\int_T f_0 \ln f_0 d\mu  -\int_T f_0 d\mu \geq \sum_{i\in I}g_i \lambda_i -\int_T e^{\sum_{i\in I}\lambda_i t^i }d\mu (t)=L(\lambda -\lambda_0 )+\sum_{i\in I}g_i \lambda_{0i}$$ where $\lambda_0 \! =\! (\lambda_{0i})_{i\in I}$ with $\lambda_{0i}\! =\! \sum_{\iota =1}^n \delta_{i,\, 2k{\rm e}_\iota }$ and
 $\delta_{i,j}$ is Kronecker's symbol.
%$\delta_{i,j}\! =\! 1$ if $i\! =\! j$ and $0$ if $i\! \not\! =\! j$. 
Since $\lambda$ was arbitrary, we get $\sup_\lambda L(\lambda )<\infty$. 
 %The existence of $\lambda^*$ such that $\sup \, L=L(\lambda^* )$ holds by general Fenchel theoretic results of dual attainment, as follows. 
 % consider the finite measure $\mu =\rho |_{T} dt$ on $T$. 
Now for the attainment $\sup L=\mbox{max}\, L$, we need Theorem \ref{maxim} as follows. Use
 $|t_j |\leq (\sum_{\iota =1}^n t_{\iota}^{2k})^{1/2k}$, %  for $1\! \leq \! j\! \leq \! n$,
%\begin{equation}
%\label{mon}
$$|t^i |=|t_{1}|^{i_1} \cdots |t_{n}|^{i_n }\leq (\sum_{\iota =1}^n t_{\iota}^{2k} +1)^{|i |/2k}\leq \sum_{\iota =1}^n t_{\iota}^{2k}+1
\mbox{ }\mbox{ }\mbox{ }\mbox{ (}|i|\leq 2k\mbox{)} 
$$ %\end{equation}
and $\nu \! +\! 1\!  \leq\!  e^\nu $  for $\nu \! =\! \sum_{\iota =1}^n t_{\iota}^{2k}$ to get $ \int_T |t^i | d\mu (t)\leq \int_T \rho dt \! <\! \infty$ for $i\! \in\! I$.
  %$\rho$  satisfies $\int_T \| t\|^{\max_{i\in I}|i|} \rho (t)dt<\infty$, % where $m=\max_{i\in I}|i|$, 
  %all monomial functions $t^i $ with $i\in I$ belong to $ L^1 (T,\rho dt)$. Then we can  
  Then let: ${\mathcal T}\! =\! T$, the measure
$\mu \! =\! \tilde{\rho} dt$, $p\! =\! \infty$, the moment functions $a_i (t)=t^i$ and  the integrand $\phi $ 
 be defined by $\phi (x)\! =\! x\ln x$ for $x>0$, $\phi (0)=0$ and $\phi (x)=+\infty$ for $x<0$. The feasibility hypotheses is fulfilled by $x=f_0$.
 The  convex conjugate  $\phi^* (y)=\sup_{x\geq 0}  (xy -x\ln x) $ of $\phi$ is given by $\phi^* (y)=e^{y-1}$ for $y\in \mathbb{R}$.
 We get the attainment $D=\sup {\mathcal L}$ for  ${\mathcal L}(\lambda )\! =\! L(\lambda \! 
-\lambda_0 ') \! +\! \sum_{i\in I}g_i \lambda_{0i}' $ where  $\lambda_0 '\! =\!  (\lambda_{0i}')_{i\in I}$ with $\lambda_{0i}'\! =\!
\lambda_{0i}\! +\! \delta_{i,0}$. Thus we  obtain a $\lambda^*$ such that $\sup L=L(\lambda^* )$. 
     
%We start with the implication (b) $\Rightarrow$ (a),  since the opposite one.
(b) $\Rightarrow$ (a)  %We prove this implication in two steps 1) and 2). 
 Let $\lambda^* \in \mathbb{R}^N$ such that $\sup L=L(\lambda^* )$. We prove that  $\varphi_g $  satisfies the positivity condition  in Theorem \ref{poz}. %, whence
 % $f$ with respect to $dt$. 
%Hence, $f:=\rho \tilde{f}$ will be a solution of (\ref{momentele}).
Let $p=\sum_{i\in I}\lambda_i X^i$, $p\not \equiv 0$ be arbitrary such that $p(t)\leq 0$ for  $t\in T$.  The vector
$\lambda :=(\lambda_i )_{i\in I}$ is then $\not =0$. For any $r> 0$, set $e_r (t)=e^{r\sum_{i\in I}\lambda_i t^i}.$ Thus $e_r (t)\leq 1$ for  $t\in T$. 
 Then the integral term $\int_T e_r \rho dt$ of  $L(r\lambda )=r\sum_{i\in I}g_i \lambda_i -\int_T e_r \rho \, dt$ remains  bounded
 as $r\to \infty$. Hence  $\varphi_g p=\sum_{i\in I}g_i \lambda_i \leq 0$, for otherwise the linear term $r\varphi_g p$
of $L(r\lambda )$ would  give $\sup L=\infty$   that is false. 
%It remains to show that the linear term is 
%actually $<0$.
Assume that $\varphi_g p =0$.
%From $\sum_{i\in I}g_i \lambda_i  =0$ and  the boundedness of the term $-\int_T e_r \rho dt $ we similarly derive 
 Then the restriction  of the  function $L$ to the  half-line $\ell :=\{ r\lambda :r>0\}$ is  given by the function $ r\mapsto -\int_T e_r \rho dt $. This function is finite, bounded and strictly
  monotonically increasing on $(0,\infty )$. Use to this aim that $0<e_r \leq 1$, $\int_T \rho dt<\infty$, $e_r =e^{rp}$ with $p\leq 0$
and  $L|_\ell $ is strictly concave.
 Then a finite limit $\lim\limits_{r\to \infty }L(r\lambda )=\sup\limits_\ell L $ exists, in particular  $\sup\limits_{r\geq 1}|L(r\lambda )|<\infty$.
%$\sum_{i\in I}g_i \lambda_i \leq 0  =0$, 
For  $a\! >\! 0$,
$$\infty >L(\lambda^* +a\lambda )=\sum_{i\in I}g_i \lambda_{i}^*  +a\sum_{i\in I}g_i \lambda_i -
\int_T e^{\sum_{i\in I}\lambda_{i}^* t^i }e^{a\sum_{i\in I}\lambda_i t^i}\rho (t)dt$$
$$\geq \sum_{i\in I}g_i \lambda_{i}^*  +r\cdot 0-
\int_T e^{\sum_{i\in I}\lambda_{i}^* t^i }\rho (t)dt=L(\lambda^* )=\max L\geq L(0)>-\infty$$
because $\sum_{i\in I}g_i \lambda_i  =0$ and $\sum_{i\in I}\lambda_i t^i \leq 0$ for all $t\in T$. Hence $L$
is finite  on every point of the half-line $\{ \lambda^* +a\lambda \}_{a> 0}$. %In particular, $L(\lambda^* +\lambda )\in \mathbb{R}$.
% Since $\max L=L(\lambda^* )$ and $L$ is strictly concave, the graph of the restriction $L|_{\sigma^*}$ 
%must remain strictly below the open half-line
% given by $H=\{ (\lambda^* +a\lambda ,L(\lambda^* )):a> 0\} \subset \mathbb{R}^N \times \mathbb{R}$.
Note that $\lambda^* $ cannot be colinear with $\lambda $ due the  behaviour  of $L$ on 
$\ell$: firstly, $\lambda^* \not \in \ell$  because  $L$ reaches 
its global maximum only in $\lambda^*$ while $L|_\ell$ increases strictly along $\ell$ as $r\to \infty$. Also $\lambda^* \not \in \{ 0\} \cup (-\ell )$, for otherwise the  concavity  of the restriction $L|_{\mathbb{R}\lambda }:\mathbb{R} \lambda \to \{ -\infty\}\cup \mathbb{R}$ of $L$ to the line $\mathbb{R} \lambda$ 
would imply, for some $r\geq 0$ with $\lambda^* =-r\lambda $, that $L(r\lambda )\geq L(0)=L(\frac{1}{2}(\lambda^* +r\lambda ))\geq \frac{1}{2} (L(\lambda^* )+L(r\lambda ))$, whence $L(\lambda^* )\leq L(r\lambda )<\sup L|_\ell  \leq \sup L=L(\lambda^* )$ that is impossible.
%; write to this aim the inequalities  
Thus $\lambda^* \not \in \mathbb{R} \lambda$.
%Set $a_0 =2$ and $r_0 =1$. 
Then a 2-dimensional drawing shows that for every $r>1$ there is a unique point $x_r$ 
of intersection of  the segments $(\lambda^* ,r\lambda )$ and $(\lambda ,\lambda^* +\lambda )$. 
%That is, $\{ x_r \} =(\lambda^* ,r\lambda )\cap (\lambda ,\lambda^* +\lambda )$
 Write to this aim $x_r =s\lambda^* +(1-s)r\lambda =s' \lambda +(1-s')(\lambda^* +\lambda )$ with coefficients $s=s_r ,\, s'=s'_r  $, 
use the linear independence of $\lambda^*$, $\lambda$ and get $s=(r-1)/r$, $s'=1-s$ whence $s,\, s'\in (0,1)$ and  $\lim_{r\to \infty}s'_r =0$. Then $\lim_{r\to \infty }x_r =\lambda^* +\lambda$.
The concavity (and hence, continuity  \cite{ro}) of $L$ on the segment $(\lambda
,\lambda^* +\lambda \, ]$ gives $\lim_{r\to \infty}L(x_r )=L(\lambda^* +\lambda )<L(\lambda^* )$ with strict inequality, because 
the point $\lambda^*$ of maximum of  $L$
is unique. 
But $L(x_r )=L(s\lambda^* +(1-s)r\lambda )\geq  sL(\lambda^* )+(1-s)L(r\lambda )$ and 
letting $r \to \infty$ we derive, using $\lim_{r\to \infty}s_r =1$ and  $\sup_{r\geq 1}|L(r\lambda )|<\infty$, that $\lim_{r\to \infty}L(x_r )\geq L(\lambda^* )$. We got a contradiction.
% with
%the strict opposite estimate from above. 
Then  $\varphi_g p <0$. The feasibility of problem (\ref{prima}) follows then by Theorem \ref{poz}.
%Then for every $r>0$, the whole closed segment $\sigma_r :=[\lambda^* ,rv] =\{ s\lambda^* +(1-s)rv\}_{0\leq s\leq 1} $
% between $\lambda^* $ and $r\lambda$ is  in the proper domain of $L$, namely $L_{\sigma_r}>-\infty$, and moreover the graph of $L_{\sigma_r}$ is above the
% chord $C_r$ relating the points $(\lambda^* ,L(\lambda^* ))$ and $(r\lambda ,L(r\lambda ))$. 
%But as $r\to \infty$  the chord $C_r$ approaches $H$ (strictly speaking, $\lim_{r \to \infty}\sup_{v\in C_r : \| v\| \leq 1}\inf_{h\in H}\| v-h\| =0$), 
%and so it follows that the graph of $L|_{\sigma^*}$ should be above $H$, that is false. Thus $\sum_{i\in I}g_i \lambda_i \leq 0  =0$.
%$\Box$
%Then obtain the dual attainment we can strengthen the conclusion of (a) and obtain also  the fact that $\sup L$ is attained. The  attainment holds  feas}
%\end{remark}
%\vspace{3 mm}
 $\Box$
\vspace{3 mm}

\noindent {\bf Remarks} 
Since $\lambda^* $ may be on the boundary of  $\mbox{dom}\, L\! :=\! \{ \! \lambda \! :\! L(\lambda )\! >\! -\! \infty\}$, 
one cannot prove (b) $\Rightarrow$ (a) by  derivating under the integral in $\lambda^*$, and the $h$-minimization  may fail \cite{junk}. 
Additional hypotheses may compel $\lambda^*$ to be interior to  $\mbox{dom}\, L$  \cite{Levermore} in which case  the entropy minimization  
  can be obtained \cite{Moreau}, providing the particular solution $f_* (t)\! =\! e^{\sum_{i\in I}\lambda_{i}^* t^i}$, see for instance \cite{MaEn}. 
For example let
 $T\! =\! \mathbb{R}^n$, $I\! =\! \{ i:|i|\leq 2k\}$ and $\rho (t)\! =\! e^{-\| t\|^{2k}}$. 
By Theorem \ref{marg},  problem (\ref{prima}) is feasible % by the existence of a maximum for $L$.
if and only if $L$ is bounded from above and attains its  maximum in a  point $\lambda^* $,
even when a minimum entropy solution does not exist. By  Fatou's lemma and Lebesgue's dominated convergence theorem, 
$f_0 \! :=\! e^{\sum_{|i|\leq 2k}\lambda_{i}^* t^i}$
  has finite moments of order $\leq 2k$, we can get  $\int t^i f_0 dt \! =\! g_i$ for $|i|\! <\! 2k$   and $\int  t_{\iota}^{2k} f_0  dt\! \leq
\!  g_{2k{\rm e}_\iota}$ ($1\leq \iota \leq n$),
 but the  equalities (\ref{prima}) may fail for $|i|=2k$ \cite{junk}. By integration in polar coordinates, the homogeneous polynomial
 $p\! :=\! \sum_{|i|=2k}\lambda_{i}^* X^i$ is shown to always satisfy $p (t)  \! \leq\! 0$ on $\mathbb{R}^n$; 
  if moreover  $p (t)  <0$ for all  $t\! \not =\! 0$, then $\lambda^*$ is interior to $\mbox{dom}\, L$ and $f_0$  is indeed  a solution of problem (\ref{prima}), 
 $f_0 =f_*$. We omit the details and refer the reader to  \cite{Levermore}, \cite{junk}. %, \cite{junk}.
  %\end{remarks}
  
Note also that whenever $\rho$ is at our disposal, various choices may be tried \cite{MaEn} to facilitate   the numerical maximization of $L=L_\rho$.
 \vspace{5 mm}
 
{\bf Acknowledgements} The present work was supported by the grants IAA100190903 of GAAV and 201/09/0473 GACR, RVO: 67985840.
   
 %, whenever condition (b) holds for a set of data

Institute of Mathematics, AS CR

Zitna 25

115 67 Prague 1

Czech Republic
\vspace{1 mm}

{\it ambrozie@math.cas.cz}
\vspace{3 mm}

and: Institute of Mathematics "Simion Stoilow" - Romanian Academy, 

\hspace{8 mm} PO Box 1-764, 014700 Bucharest, Romania

  \end{document}